%% file: main.tex
\documentclass{article}

\input{preamble/preamble}
\title{There is an equivalence relation whose von Neumann algebra is not Connes embeddable}
\author{Aareyan Manzoor \\ \texttt{a2manzoo@uwaterloo.ca}}

\begin{document}

\maketitle

 \input{abstract}

\input{prelim}

\input{IRS}

\input{Non_local_games}

\input{further}
\input{acknowledgements}

\printbibliography

\end{document}

%% file: preamble/preamble.tex
\input{preamble/Reduced_packages}

\input{preamble/operator_preamble}

%% file: preamble/Reduced_packages.tex
\usepackage{amsmath}
\usepackage{amsfonts}
\usepackage{amsthm}
\usepackage{mathtools}
\usepackage{amssymb}
\usepackage[shortlabels]{enumitem}
\usepackage[style=alphabetic-verb, maxcitenames=50,maxnames=50, backend= biber,url = false]{biblatex}
\AtEveryBibitem{%
  \clearfield{note}%
}
\usepackage[bookmarks=true, colorlinks=true, linkcolor=blue, citecolor= red,hypertexnames=false,hyperindex,breaklinks]{hyperref}
\PassOptionsToPackage{linktocpage}{hyperref}

\newtheorem{theorem}{Theorem}[section]
\newtheorem{lemma}[theorem]{Lemma}

\newtheorem{proposition}[theorem]{Proposition}
\newtheorem{definition}[theorem]{Definition}
\newtheorem{example}[theorem]{Example}

%% file: preamble/operator_preamble.tex
\renewcommand{\epsilon}{\varepsilon}

\DeclarePairedDelimiterX{\norm}[1]{\lVert}{\rVert}{#1}
\DeclarePairedDelimiterX{\abs}[1]{\lvert}{\rvert}{#1}
\DeclarePairedDelimiter\bra{\langle}{\rvert}
\DeclarePairedDelimiter\ket{\lvert}{\rangle}
\DeclarePairedDelimiterX\braket[2]{\langle}{\rangle}{#1\,\delimsize\vert\,\mathopen{}#2}
\newcommand{\ip}[3]{\bra{#1}#2\ket{#3}}
\DeclarePairedDelimiterX{\gen}[1]{\langle}{\rangle}{#1}
\DeclareMathOperator{\tr}{tr}

\DeclareMathOperator{\Val}{Val}
\newcommand{\sub}{\mathfrak{sub}}
\newcommand{\IRS}{\operatorname{IRS}}
\DeclareMathOperator{\Stab}{Stab}
\DeclareMathOperator{\Aut}{Aut}
\usepackage{graphicx}

\newcommand{\Conv}{\mathop{\scalebox{1.5}{\raisebox{-0.2ex}{$\ast$}}}}

%% file: abstract.tex
\begin{abstract}
          The landmark quantum complexity result MIP$^*$=RE was used to prove the existence of a non Connes embeddable tracial von Neumann algebra.  Recently, similar ideas were used to give a negative solution to the Aldous-Lyons conjecture: there is a non co-sofic IRS on any non-abelian free group. We define a notion of hyperlinearity for an IRS and show that there is a non co-hyperlinear IRS on any non-abelian free group.  As a corollary, we prove that there is a relation whose von Neumann algebra is not Connes embeddable. We do this by significantly simplifying the reduction of Aldous-Lyons to non-local games, removing the need for subgroup tests entirely.
 
\end{abstract}

%% file: prelim.tex
\section{Introduction}
One of the biggest problems in von Neumann algebras was the Connes embedding problem (CEP), which asked if each tracial von Neumann algebra embedded into an ultrapower $R^U$ of the hyperfinite II$_1$ factor $R$ in a trace preserving way. This is a finite dimensional approximation property. This problem was famous because it turned out to be equivalent to several different problems all over operator algebras, but also in quantum information theory. See \cite{goldbring_connes_2021} for detailed discussions on different equivalent formulations of this problem. 

The landmark quantum complexity result MIP*=RE \cite{ji_mipre_2022} showed a negative solution to the CEP through a long chain of equivalences. This proof was very non-constructive, essentially showing that if each tracial von Neumann algebra was Connes embeddable: then the Halting problem would be decidable. The idea is that CEP is equivalent to a separation of convex sets problem (see Proposition \ref{prop:2.15}), and certain nice convex sets can be separated using a nice computability trick involving linear functionals (see the proof of Theorem \ref{thm:3.15}). As this proof is non-constructive, it can be hard to say anything more about about non-Connes embeddable factors, as there is not a nice counterexample to test properties on.

However it has been possible to use the complexity result MIP*=RE to build non-Connes embeddible II$_1$ factors with additional structure. For example, Goldbring and Hart \cite{GoldbringHart2024UniversalTheoryR} use MIP*=RE to prove the existence of non-Connes embeddable II$_1$ factors satisfying certain classes of model theoretically nice properties. We will use the recent strengtening of MIP*=RE, the so-called TailoredMIP*=RE \cite{bowen_aldous--lyons_2024}, to get the main result in Theorem \ref{thm:1.1}. TailoredMIP*=RE was also used to prove that both satisfying L\"uck's determinant conjecture \cite{man25} and Gottshack's surjunctivity conjecture \cite{bc25} are distinct from soficity for invariant random subgroups.

To each group $\Gamma$ one can associate a tracial von Neumann algebra $L(\Gamma)$, and one question is whether the class of group von Neumann algebras has a non-Connes embeddable algebra. The groups whose von Neumann algebra is Connes embeddable are called hyperlinear.
Another finite approximation property for groups is soficity. Roughly speaking: hyperlinearity gives the group an approximation by finite dimensional matrix groups, while soficity gives an approximation by finite groups. Each sofic group is in particular hyperlinear \cite{elek_hyperlinearity_2014}. So an easier problem than asking ``Is there a non-hyperlinear group?'' is ``Is there a non-sofic group?''. While these computability techniques can't solve these problems yet on the level of groups, they can solve these at the level of group actions:

\begin{theorem}\label{thm:1.1}
    There is some Property $(T)$ ergodic countable probability measure preserving Borel equivalence relation $\mathcal{R}$ so that the property $(T)$ II$_1$ factor $L(\mathcal{R})$ is not Connes embeddable.
\end{theorem}
To understand how TailoredMIP*=RE can be used to solve this, we first recall how MIP*=RE was used to solve the CEP. 
\subsection{Separation of convex sets}

Fix $\Gamma=\mathbb F_{n,2^m}:= \langle u_1\dots u_n: u_i^{2^m}=1\rangle$ be the group generated by $n$ free elements of order $2^m$. Let $T(\Gamma)$ denote the compact convex set of traces(i.e positive definite conjugation invariant functions) on $\Gamma$, and let $AT(\Gamma)\subseteq T(\Gamma)$ be the subset of amenable traces. By work of Kirchberg, if $AT(\Gamma)\neq T(\Gamma)$ for some $(n,m)$, then there exists a non-Connes embeddable $\mathrm{II}_1$ factor. See Proposition \ref{prop:2.15}.

We will witness the strict inclusion via functionals corresponding to non-local games. It is not important for this sketch to know what non-local games are, only that their strategies come from traces on some $\Gamma$ (The details are discussed in Section \ref{sec:3.1}). So for each non-local game $\mathfrak{G}$ with question set of size $n$ and answer set of size $m$, and each trace $\tau$ on $\Gamma$, we can define $\Val(\mathfrak{G},\tau)\in [0,1]$ as the winning probability under the strategy corresponding to $\tau$.
Define
\[
\omega^{co}(\mathfrak G)=\sup_{\tau\in T(\Gamma)}\Val(\mathfrak G,\tau),
\qquad
\omega^{*}(\mathfrak G)=\sup_{\tau\in AT(\Gamma)}\Val(\mathfrak G,\tau).
\]
Thus a strict gap $\omega^{co}(\mathfrak G)>\omega^{*}(\mathfrak G)$ produces a trace in $T(\Gamma)\setminus AT(\Gamma)$, and hence a non-Connes embeddable factor.

Computability enters because both quantities admit effective one-sided approximations: there are computably enumerable lower bounds $\alpha_d(\mathfrak G)\nearrow \omega^{*}(\mathfrak G)$ and computably enumerable upper bounds $\beta_d(\mathfrak G)\searrow \omega^{co}(\mathfrak G)$. The theorem MIP$^*=$RE \cite{ji_mipre_2022} implies that $\omega^{*}(\mathfrak G)$ is not computable as a function of $\mathfrak G$.

If one had $\omega^{co}(\mathfrak G)=\omega^{*}(\mathfrak G)$ for all games, then running the procedures for $\alpha_d$ and $\beta_d$ in parallel and waiting until $\beta_d(\mathfrak G)-\alpha_d(\mathfrak G)<2^{-k}$ would yield an algorithm to approximate $\omega^{*}(\mathfrak G)$ to precision $2^{-k}$, contradicting uncomputability. Therefore there exists some game $\mathfrak G$ with $\omega^{co}(\mathfrak G)>\omega^{*}(\mathfrak G)$, and hence some $(n,m)$ with $AT(\Gamma)\neq T(\Gamma)$.

\subsection{Invariant Random subgroups}
An invariant random subgroup (IRS) of a group $\Gamma$ is a Borel conjugation invariant probability measure on its set of subgroups. These objects generalize the notion of a normal subgroup, as the Dirac measure of a normal subgroup is an IRS. Since normal subgroups of the free groups correspond to all groups by taking quotients, IRS's on a free group generalize the notion of a group. The notion of a co-sofic IRS is the generalization of a sofic group in this framework, see Proposition \ref{prop:1.14}.

The existence of a non co-sofic IRS on a free group was proven in \cite{bowen_aldous--lyons_2024-1},\cite{bowen_aldous--lyons_2024}, disproving the Aldous-Lyons conjecture. To do this, they first in \cite{bowen_aldous--lyons_2024-1} define subgroup tests, a type of combinatorial game that are very natural probability theoretically. IRS's on free groups correspond to strategies on these, and by the same technique discussed for CEP, they reduced this problem to a computability problem. Then they reduced that computability problem to proving a stronger version of MIP*=RE, and in the second paper proved this computability result TailoredMIP*=RE \cite{bowen_aldous--lyons_2024}.

The reduction from subgroup tests in \cite{bowen_aldous--lyons_2024-1} to non-local games is very long and technical. Fortunately, our approach avoids this entirely, and directly disproves Aldous-Lyons by non-local games. Taking \cite{bowen_aldous--lyons_2024} as a blackbox, we are able to reprove the main result of \cite{bowen_aldous--lyons_2024-1}: we only use some minor results from there to prove the computable upper bounds for $\omega_{IRS}(\mathfrak{G})$ in Section \ref{sec:3.2}, otherwise this paper is independent of \cite{bowen_aldous--lyons_2024-1}.

 The idea is that each IRS on a free group naturally corresponds to a trace on it:
\[\tau_H(g) = \mathbb{P}(g\in H).\] This then naturally corresponds to strategies on non-local games. This allows us to define $\omega_{IRS}(\mathfrak{G})$ to be the best winning strategy over traces coming from IRS's (with a minor technicality). This directly converts Aldous-Lyons to a problem about non-local games, without the need to go through subgroup tests of \cite{bowen_aldous--lyons_2024-1}.

Our approach also allows us to upgrade soficity to hyperlinearity. we define for each IRS $H$ on a group $\Gamma$, a von Neumann algebra $L(\Gamma/H)$ generalizing the notion of quotienting against a normal subgroup. It is then natural to call an IRS $H$ of a free group co-hyperlinear if $L(\Gamma/H)$ is Connes embeddable. We prove:
\begin{theorem}\label{thm:1.0}
    There is some free group with a non co-hyperlinear invariant random subgroup.
\end{theorem}

In particular, each $L(\Gamma/H)$ is a subalgebra of $L(\mathcal{R})$(see discussion after Definition \ref{def:1.14}) for some countable Borel equivalence relation $\mathcal{R}$. This will allow us to recover Theorem \ref{thm:1.1} from Theorem \ref{thm:1.0}

%% file: IRS.tex
\section{Invariant Random Subgroups}\label{sec:2}

This subsection will set up the terminology we use for invariant random subgroups. Most of these are standard definitions, compare with \cite{bowen_aldous--lyons_2024-1}.
\begin{definition}
    Let $\Gamma$ be a discrete group. Let $\sub(\Gamma)$ denote the space of subgroups of $\Gamma$ topologized as a subset of $\{0,1\}^\Gamma$. Define a \textbf{Random subgroup} as a Borel probability measure on $\sub(\Gamma)$. Further, if the measure is invariant under the conjugation action of $\Gamma$ on $\sub(\Gamma)$ then we call it an \textbf{Invariant Random Subgroup(IRS)}. We will denote by $\IRS(\Gamma)$ the convex set of IRS's on $\Gamma$.
\end{definition}
Note that $\sub(\Gamma)$ is compact, so we can talk about the weak$^*$ topology on $\IRS(\Gamma)$.  The topology on $\sub(\Gamma)$ can be understood as follows: a net of subgroups $H_\lambda$ converges to $H$ means 
\[g\in H \iff g\in H_\lambda \text{  eventually}.\]
The topology on $\IRS(\Gamma)$ can be understood as follows: a net $\mu_\lambda$ converges to $\mu$ if for each continous $f: \sub(\Gamma) \to \mathbb{R}$, we have
\[\int f d\mu_\lambda \to \int f d\mu.\]
Note in particular this is compact by the Banach-Alaoglu theorem.

Here are some examples of IRS's:
\begin{example}
    Let $\Gamma$ be a discrete group. Let $N\subset \Gamma$ be a normal subgroup, then the point mass $\delta_N$ is an IRS.
\end{example}
So we should think of an IRS as a generalization of a normal subgroup of $\Gamma$. We will only really care about IRS's of free groups, and since normal subgroups of free groups correspond to all finitely generated groups by looking at the quotients: we will consider these IRS as a generalization of the notion of a group by some generalized notion of a quotient (see Definition \ref{def:1.11}).
\begin{example}    
 Let $(X,\mu)$ be a probability space, and $\alpha:\Gamma\to \operatorname{Aut}(X,\mu)$ be a probability measure preserving (p.m.p.) action. Then we can describe an IRS as follows: first sample a $x\in X$ according to $\mu$, and then look at \[\Stab(x)= \{g\in \Gamma: gx=x\}.\] We will call the induced probability measure on $\sub(\Gamma)$ by  $\Stab(\alpha)$. Formally if we look at the map $\Stab: X\to \sub(\Gamma)$ sending a point to its stabilizer, $\Stab(\alpha)$ is the pushforward of $\mu$ under this map. 
\end{example}
Note that the action $\alpha$ being probability measure preserving gives conjugation invariance of the $\Stab(\alpha)$. 
%\mbnote{Briefly explain why this is a conjugation invariant measure.}

%\mbnote{Why do we need $X$ compact?  Is the terminology compact Borel space standard?  Are you implicitly assuming $\mu(X) = 1$?  It seems you are not making this assumption in the proof of Lemma 1.21}
A remarkable result is that all IRS's are of this form:
\begin{theorem}\cite[Proposition 14]{abert_kestens_2014}
    Let $\Gamma$ be a countable discrete group, and $H\in \IRS(\Gamma)$. Then there is a p.m.p. action $\alpha$ of $\Gamma$ on some compact space $X$ with Borel probability measure $\mu$ such that $H=\Stab(\alpha)$. If $H$ was ergodic, then $\alpha$ can be picked to be ergodic.
\end{theorem}
We can now define:%\mbnote{add refs to where these defns come from}:
\begin{definition}
    Let $\Gamma$ be a discrete group, and $\alpha$ a p.m.p. action of it on a finite set. Then the IRS $\Stab(\alpha)$ will be called \textbf{finitely described}. A weak$^*$ limit of finitely described IRS's will be called co-sofic.
\end{definition}
The name is motivated by:
%\mbnote{add ref}:
\begin{proposition}\label{prop:1.14}\cite[lemma 16]{abert_rank_2017}
    Let $\Gamma$ be a free group, and $N$ a normal subgroup. Then $\delta_N$ is co-sofic as an IRS $\iff$ $\Gamma/N$ is sofic as a group.
\end{proposition}
So if we think of IRS's of free groups as an extension of the notion of a group, we should think of co-sofic IRS's of the free group as an extension of the notion of a sofic group. One of the biggest open problems in this area is: Is every group sofic? What was proven in \cite{bowen_aldous--lyons_2024-1} is:
\begin{theorem}
    Let $\Gamma$ be a non-abelian free group. Then there is a non co-sofic IRS on it.
\end{theorem}
%\mbnote{By countable, you mean countable generating set?  Countably many elements? For the latter you'd have to exclude $\mathbb Z$.  Need to clarify.}
This of course does not prove the existence of non-sofic groups: indeed this non co-sofic IRS doesn't have to be a point-mass. 

A slightly harder problem than finding a non-sofic group is finding a non hyperlinear group. The rest of the paper is dedicated to setting up what co-hyperlinearity should mean for IRS's and proving there is a non co-hyperlinear IRS on a free group.
%\mbnote{Again mixing co- and non-co-...}

\subsection{The algebra of an IRS}
We will formalize the notion of quotienting against an IRS in this section. The idea is to go through traces on the universal group $C^*$ algebra $C^*(\Gamma)$ of $\Gamma$. This is a $C^*$-algebra densely containing the group ring $\mathbb{C}[\Gamma]$ that encodes all the representation theoretic date of $\Gamma$. More precisely, there is a correspondence between representations of $C^*(\Gamma)$ on Hilbert spaces and unitary representations of $\Gamma$ on the same Hilbert space. See \cite[Section 2.5]{brown_c-algebras_2008} for exposition.
\begin{definition}
    Let $\Gamma$ be a discrete group, and $H\in \IRS(\Gamma)$. Then one can define a trace $\tau_H$ on $C^*(\Gamma)$ by:
    \[\tau_H(g) := \mathbb{P}(g\in H).\]
\end{definition}
%\mbnote{You should say somewhere what the notation  $C^*(\Gamma)$ stands for.}
Note that invariance gives $\tau_H$ conjugation invariance, but one still needs to check positivity. We get this from Proposition \ref{prop:1.13}. If $H=\Stab(\alpha)$ where $\alpha$ is a p.m.p. action of $\Gamma$ on some Borel probability space $(X,\mu)$, then 
\[\tau_H(g) = \mu(x: gx=x).\]

In group theory these objects are called  characters of $\Gamma$ (positive definite conjugation invariant functions). Naimark's dialation theorem gives this equivalence, see \cite[theorem 4.8]{paulsen_completely_2003} for exposition.

Denote the trace complex of $\Gamma$ by $T(\Gamma)$, and equip it with the weak$^*$ topology. This is the topology of pointwise convergence on $C^*(\Gamma)$, which is the same as pointwise convergence when restricted to $\Gamma$. Then note there is an affine map $\IRS(\Gamma)\to T(\Gamma)$ given by $\IRS(\Gamma) \owns H \mapsto \tau_H \in T(\Gamma)$. This map is neither injective nor surjective generally.

Not being surjective is clear as a trace $\tau \in T(\Gamma)$  can be negative on a group element, which is not the case for $\tau_H$ for $H \in \IRS(\Gamma)$. For non-injectivity consider $S_3$, which has $4$ conjugacy classes of subgroups but $T(S_3)$ has $3$ irreducible characters. An IRS on a finite group is simply a probability measure on the set of conjugacy class of subgroups, while a character on a finite group is in the convex hull of the irreducible characters. So a dimension arguement forces this map to not be injective. This map is however $w^*-w^*$ continuous, as $\{H: g\in H\}$ is a clopen set making its indicator continuous. So the image of this, call it $T_{IRS}(\Gamma)$, is compact.

One of the main themes in this paper will be looking at the von Neumann algebra that corresponds to a trace. This is done by looking at the Gelfand-Naimark-Segal(GNS) representation of a trace.

Let $\mathcal{A}$ be a $C^*$-algebra and $\tau$ a trace on it, then the GNS representation of $\tau$ is a representation $\pi_\tau: \mathcal{A}\to B(H_\tau)$ and a vector $\xi \in H_\tau$ with $\tau(a) = \ip{\xi}{\pi(a)}{\xi}$. The von Neumann algebra of $\tau$ will be $\pi_\tau(\mathcal{A})'' \subset B(H_\tau)$, and note that $\tau$ extends to this algebra by
\[\tilde{\tau}(a) =  \ip{\xi}{\pi(a)}{\xi}\quad a\in \pi_\tau(\mathcal{A})''.\]
 This algebra is well defined due to uniqueness of GNS: any two cyclic representation inducing $\tau$ will be unitarily equivalent. We will denote this von Neumann algebra as $\tau(\mathcal{A})''$. If $\mathcal{A}=C^*(\Gamma)$, we will denote this as $\tau(\Gamma)''$. It will be implicit that we are considering this von Neumann algebra with the extension of $\tau$ as its trace.

  Consider the left regular representation $\lambda:\Gamma \to U( \ell^2(\Gamma))$, where $(\lambda(g)f)(h) = f(g^{-1}h)$ for $g,h \in \Gamma$. Then we define the group von Neumann algebra to be $\lambda(\Gamma)''\subset B(\ell^2(\Gamma))$. We consider it with the trace $\tau(g) = 1$ if $g=e$ and $0$ otherwise on group elements.   
%\mbnote{This is a very long definitions with a lot of terminology packed in.  I would just convert this to a paragraph rather than a definition.}

Before defining our IRS algebras, let us make one observation:
\begin{proposition}
    Let $\Gamma$ be a discrete group and $N$ a normal subgroup. If $\tau_N$ is the trace on $\Gamma$ corresponding to the point mass $\delta_N \in \IRS(\Gamma)$, then 
    \[\tau_N(\Gamma)''= L(\Gamma/N)\]
    where $L(\Gamma/N)$ is the group von Neumann algebra of $\Gamma/N$ with its standard trace.
\end{proposition}
%\mbnote{Did you define standard trace anywhere? Did $L(G)$ get defined?}
\begin{proof}
    First note that $\tau_N(g) = 1$ if $g\in N$ and $0$ otherwise.  Now take the representation:
    \[\pi: \Gamma \to U(\ell^2(\Gamma/N))\]
    where $(\pi(g)f)([h]) = f([g^{-1}h])$ for $g\in \Gamma, [h]\in \Gamma/N$. Note that all elements of $N$ act trivially, so this is actually a homomorphism $\tilde{\pi}: \Gamma/N \to U(\ell^2(\Gamma/N))$. Note this is just the left regular representation of $\Gamma/N$.
    
    Fix the vector $\xi \in \ell^2(\Gamma/N)$ that is $1$ on the identity coset and $0$ everywhere else. Then it is clear
    \[\tau_N(g) = \ip{\xi}{\pi(g)}{\xi}\]
    So the corresponding algebra is $\pi(\Gamma)'' = \tilde{\pi}(\Gamma/N)''$ and by definition this is the group von Neumann algebra $L(\Gamma/N)$. Note that the extension of the trace is $1$ on elements of $N$ and $0$ on every other group element, i.e the standard trace on the group von Neumann algebra.
\end{proof}
Based on this, we can now define the quotient of a group by an IRS:
\begin{definition}\label{def:1.11}
    Let $\Gamma$ be a discrete group and $H\in \IRS(\Gamma)$. Then define the \textbf{quotient of $\Gamma$ by $H$} as the tracial von Neumann algebra:
    \[L(\Gamma/H):= \tau_H(\Gamma)''.\]
\end{definition}
%\mbnote{Is it simply ``the quotient'' or ``the quotient of $\Gamma$ by the IRS $H$''?  Also, I noticed some inconsistency in your use of boldface/nonboldface for definitions.  Make sure to follow one convention.}
A natural question now is this: can we write the GNS representation of $\tau_H$ explicitly? The answer is yes. 

To see this, let $H\in \IRS(\Gamma)$ correspond to the p.m.p. action $\alpha$ of $\Gamma$ on some $(x,\mu)$. Let $\mathcal{R}_\alpha\subset X\times X$ be the relation defined by the orbits, i.e $x\equiv gx$. Let $\pi_0: \mathcal{R}_\alpha\to X$ be given by $(x,y)\mapsto x$ and let $\Gamma$ act on $\mathcal{R}_\alpha$ on the first coordinate, i.e $g(x,y)= (gx,y)$. Define 
 \[\nu(E) := \int |\pi_0^{-1}(x) \cap E| d \mu(x).\] This is a $\sigma$-finite measure on $\mathcal{R}_\alpha$ \cite[section 1.5]{anantharaman_introduction_nodate}.
    %\mbnote{I don't think $\nu$ needs to be given the status of a definition - I'd turn it into a paragraph preceding the next proposition.  You should also define the representation $\pi$ prior to the next result.}
 \begin{proposition}\cite[Theorem 9]{vershik_nonfree_2010}\label{prop:1.13}
        Let $H$ be an IRS of $\Gamma$. Suppose it corresponds to p.m.p. action $\alpha$ of $\Gamma$ on $(X,\mu)$,i.e $H=\Stab(\alpha)$. Let $\pi: \Gamma\to U(L^2(\mathcal{R}_\alpha,\nu))$ correspond to the action described above. Let $\Delta\subset \mathcal{R}_\alpha$ be the diagonal of $X$. Then
        \[\tau_H(g) = \ip{1_\Delta}{\pi(g)}{1_\Delta}.\]
        In particular, the algebra $L(\Gamma/H) = \pi(\Gamma)''$.
    \end{proposition}
    \begin{proof}
        Note that $\pi(g)1_\Delta$ is the characteristic function on $g^{-1}\Delta=\{(g^{-1}x,x):x\in X\}$. So
        \[\begin{aligned}
        \ip{1_\Delta}{\pi(g)}{1_\Delta} &= \int_{\mathcal{R}} 1_{\Delta}1_{g^{-1}\Delta} d\nu\\
        &= \int_{\mathcal{R}} 1_{\{(x,x): x=gx\}} d\nu\\
        &= \nu\{(x,x): x=gx\}\\
        &= \int_X |\pi_0^{-1}(x) \cap \{(x,x): x=gx\}| d\mu\\
        &= \int_X 1_{x=gx} d\mu\\
        &= \mu(x: gx=x).
        \end{aligned}\]
    \end{proof}
    This can be extended to a larger algebra:
    \begin{definition}\label{def:1.14}
        Let $\alpha$ be a p.m.p. action of $\Gamma$ on $(X,\mu)$, and let $(\mathcal{R}_\alpha,\nu)$ be as in proposition \ref{prop:1.13}. For $f\in L^\infty(X,\mu)$ define $M_f \in B(L^2(\mathcal{R}_\alpha,\nu))$ as
        \[(M_f \xi) (x,y) = f(x) \xi(x,y).\]
        Then the von Neumann algebra generated by $\pi(\Gamma)$ and $\{M_f: f\in L^\infty(X,\mu)\}$ is defined to be von Neumann algebra $L(\mathcal{R}_\alpha)$ of the relation $\mathcal{R}_\alpha$. We consider this with the trace
        \[\tau_\alpha(a) = \ip{1_\Delta}{a}{1_\Delta},\quad a\in L(\mathcal{R}_\alpha)\]
    \end{definition}
    This can be defined for any measurable relation on a Borel space, but essentially reduces to this case. These algebras are useful as the inclusion $L^\infty(X)\subset L(\mathcal{R}_\alpha)$ is a complete invariant for the action $\alpha$ (up to orbit equivalence). These are also canonical examples of Cartan Inclusions.  Note in particular $L(\Gamma/\Stab(\alpha))\subset L(\mathcal{R}_\alpha)$. Also note that if $\alpha$ is ergodic then $L(\mathcal{R}_\alpha)$ is a $II_1$ factor. See \cite[section 1.5, 12.1]{anantharaman_introduction_nodate} for a discussion of all of these.

    \subsection{Co-hyperlinear IRS}

    Now that we can talk about IRS's algebraicly, it is easy to define co-hyperlinearity for it. Recall that a group $\Gamma$ is called hyperlinear if its von Neumann algebra $L(\Gamma)$ is Connes embeddable, i.e there is a trace preserving embedding $L(\Gamma)\hookrightarrow R^\omega$ of the group algebra into an ultrapower of the hyperfinite $II_1$ factor. So we want co-hyperlinearity of a $H\in \IRS(\Gamma)$ to line up with Connes embeddebility of $L(\Gamma/H)$.
    
    A crucial point for us is that a tracial von Neumann algebra being Connes embeddable is something we can express in terms of traces:
     \begin{definition}\label{def:1.16}
        A trace $\tau$ on a C* algebra $\mathcal{A}$ is called amenable if there are contractively completely positive (c.c.p) maps $\phi_n: \mathcal{A}\to M_{k_n}(\mathbb{C})$ with
        \[\norm{\phi_n(ab)-\phi_n(a)\phi_n(b)}_2 \to 0 \quad  \forall a,b\in \mathcal{A}\]
        so that $\tau(a) = \lim \tr_{k_n} \circ \phi_n(a)$. Here $\norm{A}_2^2 = \tr_n(A^*A)$ is the normalized Hilbert-Schmidt norm.
    \end{definition}
    %\mbnote{Define $\|\cdot\|_2$. Add words between defns and propositions to make the paper flow better.}
    This is a finite dimensional approximation property for the trace, and it turns out the corresponding finite dimensional approximation property for the algebra is Connes embeddability. We state the following proposition for groups satisfying either the local lifting property or Hilbert Schmidt stability. The important fact is that free groups have both of these properties \cite[Theorem 13.1.3]{brown_c-algebras_2008},\cite[Theorem 1]{HadwinShulman2018JFA}.
    \begin{proposition}
        Let $\Gamma$ be a group such that $C^*(\Gamma)$ has the Local Lifting Property, and $\tau$ a trace on it. The following are equivalent by \cite[Proposition 6.3.4]{brown_invariant_2006}:
        \begin{enumerate}[1.]
            \item $\tau(\Gamma)''$ has a trace preserving embedding into $R^\omega$
            \item $\tau$ is amenable.
        \end{enumerate}
        In addition, if $\Gamma$ is Hilbert-Schmidt stable then the following the above are also equivalent to the following by \cite[Theorem 3]{HadwinShulman2018JFA}:
        \begin{enumerate}[1.]\setcounter{enumi}{2}
            \item There is a sequence of finite dimensional representations $\pi_n:\Gamma \to U_n(\mathbb{C})$ with $\tau$ being the weak$^*$ limit of $\operatorname{tr}_n\circ \pi_n$.
            \item There is a sequence of traces $\tau_n$ with finite dimensional GNS representation  whose weak$^*$ limit is $\tau$.
        \end{enumerate}
        Here $R^\omega$ is an ultrapower of the hyperfinite $II_1$ factor.
    \end{proposition}
    Note the above proposition is not true for a general $C^*$-algebra $\mathcal{A}$. In general amenability of the trace means $\tau(A)''$ has an embedding into $R^\omega$ with a completely positive lift onto $\prod_\mathbb{N} R$. Local lifting property gives this lift for free. See \cite[section 13.1]{brown_c-algebras_2008} for a detailed discussion of the local lifting property. We will use this property in our proof of Theorem \ref{thm:3.5} in a minor way. Hilbert Schmidt stability allows the completely positive lift of the embedding of $\tau(A)''\hookrightarrow R^\omega$ to actually be a $*$-homomorphism. This gives the equivalence of 1 and 3. 
    
    Let us call the collection of traces  with finite dimensional GNS represention the finite dimensional traces of $\Gamma$, and denote this set by $T_{f.d}(\Gamma) \subset T(\Gamma)$.

    We also remark that while the traces are called amenable, the corresponding algebras are not amenable von Neumann algebras. For that, one would need to restrict to a class of traces called uniformly amenable traces. 
    
    Of course one of the nice things about the free group is that it is an universal object. It turns out every tracial von Neumann algebra arises from a trace on some free group (the extreme points correspond to factors). So the fact that there is a non-Connes embeddable tracial von Neumann algebra can be stated as:
    \begin{proposition}\label{prop:2.15}
        There is a non-amenable trace on each non-abelian free group. That is if $\Gamma$ is free then $T(\Gamma) \neq \overline{T_{f.d.}(\Gamma)}^{w^*}$.
    \end{proposition}
    See \cite[Proposition 6.3.5]{brown_invariant_2006} for more details. Now it is clear how to generalize the notion of hyperlinearity to IRS's:
    \begin{definition}
        Let $\Gamma$ be a discrete group, and $H\in \IRS(\Gamma)$. Then call $H$ \textbf{co-hyperlinear} if the corresponding trace $\tau_H$ is amenable.
    \end{definition}
    \begin{proposition}\label{prop:2.17}
        Let $\Gamma$ be a free group, and $H\in \IRS(\Gamma)$. Then
        \begin{enumerate}[1.]
            \item $H$ is co-hyperlinear $\iff$ $L(\Gamma/H)$ is Connes embeddable
            \item $H=\delta_N$ is co-hyperlinear $\iff$ $\Gamma/N$ is hyperlinear as a group.
        \end{enumerate}
    \end{proposition}
    Let us make explicit what the difference between co-hyperlinear and co-sofic is:
    \begin{lemma}\label{lemma:2.21}
        Let $\Gamma$ be a discrete group and $H\in \IRS(\Gamma)$ be ergodic. 
        
        Then $H$ being finitely described implies $\tau_H$ is a finite dimensional trace.

        Further if $\Gamma$ is finitely generated, then $\tau_H$ being a finite dimensional trace implies $H$ is finitely described.
    \end{lemma}
    \begin{proof}
        $(\implies)$ Let $H=\Stab(\alpha)$ where $\alpha$ is the action of $\Gamma$ on a finite set with the uniform measure. Then the orbit relation $\mathcal{R}_\alpha\subset X\times X$ is also a finite set and hence by proposition \ref{prop:1.13} there is a finite dimensional representation $\pi:\Gamma \to U(L^2(\mathcal{R}_\alpha))$ inducing $\tau_H$.

        $(\impliedby)$ We can find an ergodic action $\alpha$ of $\Gamma$ on some probability space $(X,\mu)$ with $H=\Stab(\alpha)$.

        Let $F_n\subset X$ denote the set of points with orbits of size less than $n$. Note that a point having orbit of size $n$ is the same as the stabilizer $K$ having index $n$. Since the action $\Gamma$ on $\Gamma/K$ uniquely determines $K$ (by looking at the stabilizer of the identity coset), the corresponding homomorphism $\Gamma \to \operatorname{Sym}(n)$ uniquely determines $K$. In particular, since $\Gamma$ is finitely generated, there are finitely many subgroups of index $n$ of $\Gamma$. So:
        \[F_n = \bigcup_{K: |\Gamma/K|\leq n} \{x : gx=x \forall g\in K\}.\]
        This is a finite union of measurable sets, and hence measurable. $F_n$ is $\Gamma$ invariant and hence has measure $0$ or $1$. If it has measure $1$, then $H$ is supported on subgroups of index $\leq n$, and this is a finite set, call it $\operatorname{Supp}(H)$. Now ergodicity of $H$ forces $\operatorname{Supp}(H)$ to be the set of conjugates of a given subgroup $K'$, and this is the IRS induced by the action of $\Gamma$ on the finite set of cosets $\Gamma/K'$. 

        If all $F_n$ have measure $0$, then the action has infinite orbits almost everywhere. By \cite[Proposition 3.1]{peterson_character_2014}, the representation $\Gamma \to U(L^2(\mathcal{R}_\alpha,\nu))$ is weakly mixing, i.e it has no finite dimensional sub-representations. In particular, $\tau_H$ cannot be finite dimensional, ruling this case out.
    \end{proof}
    %Note that if we dropped the ergodicity condition, then we could use this to show that if an IRS trace is a norm limit of finite dimensional traces (a uniformly amenable trace), it corresponds to a co-sofic IRS. The GNS of these uniformly amenable traces are hyperfinite, which is the von Neumann algebra version of an amenable group. So this would be saying ``co-amenable'' IRS are co-sofic. Amenable groups were the first class of groups proved to sofic.

    This immediately gives us:
    \begin{proposition}\label{prop:2.22}
    Let $\Gamma$ be a countable discrete group and $H \in \IRS(\Gamma)$ be co-sofic. Then $H$ is co-hyperlinear
    \end{proposition}
    \begin{proof}
        Let $H$ be co-sofic, then there is a sequence of finitely described IRS $(H_n)$ with $H_n\to H$. Since the map sending an IRS to its trace is $w^*-w^*$ continuous, we have $\tau_{H_n}\to \tau_H$. By Lemma \ref{lemma:2.21} we get that each $\tau_{H_n}$ is finite dimensional, and so $\tau_H$ is a $w^*$ limit of finite dimensional traces. By Definition \ref{def:1.16}, since these finite dimensional traces factor through matrix algebras by $*$-homomorphisms (which are in particular c.c.p maps that are approximately multiplicative), we have $\tau_H$ is amenable.
    \end{proof}
    
    Let $\Gamma$ be a free group. Note that showing a seperation between co-sofic IRS's and all IRS's is implied by showing a seperation between $T_{IRS}(\Gamma)$ and $\overline{T_{IRS}(\Gamma)\cap T_{f.d.}(\Gamma)}$. On the other hand, showing seperation between co-hyperlinear IRS's and all IRS's is showing a seperation between $T_{IRS}(\Gamma)$ and $T_{IRS}(\Gamma)\cap \overline{T_{f.d.}(\Gamma)}$. The latter is an a priori harder problem as one might be able to approximate a IRS trace with finite dimensional traces, but not finite dimensional traces arising from IRS's. Actually these notions agreeing would imply hyperlinearity and soficity agree for groups, by specializing to dirac measures.

    We also remark that the set of co-hyperlinear IRS's of free group $\Gamma$ is a $w^*$ closed convex subset of $\IRS(\Gamma)$ containing the co-sofic IRS's. This is because it is the pre-image of $T_{IRS}(\Gamma)\cap \overline{T_{f.d.}(\Gamma)}$ under the affine continuous $H\mapsto \tau_H$.

%% file: Non_local_games.tex
\section{Non-Local Games}\label{sec:3.1}

%\mbnote{Fluff it up with some words explaining what's gonna happen in this section.  Make sure to add a standard reference for the non expert reader to look things up (there's a recent preprint sruvery by Harris-Paulsen that could be mentioned).  Also, are you using the A-L paper notation here.  You should mention whose notational conventions you are using.} \\ \\
To each trace on a certain group $C^*$ algebra to be defined below, one can associate a strategy on non-local games. Since we have associated to each IRS a trace, we can associate to each IRS a strategy for some non-local game. See \cite[section 4,5]{goldbring_connes_2021} for a less bare bones discussion of the topic. 
%\mbnote{$X = Q$ throughout this section.  Make sure to keep the notation consistent.}
\begin{definition}
    A non-local game $\mathfrak{G}$ has the following parameters: a finite question set $Q$, a probability distribution $q$ on $Q\times Q$, an integer $m$ describing the length of the answer, and a decider function $D:(\{0,1\}^{m})^2\times Q^2 \to \{0,1\}$. We will write this as $D(a,b|x,y)$ for $a,b \in \{0,1\}^m$ and $x,y \in Q$.

    We will also associate to this game a set 
    \[S(\mathfrak{G})=\{u_{x,i}: x\in Q, 1\leq i\leq m\}\]
    each $u_{x,i}$ is to be regarded as a formal variable which encodes the value of the $i$th bit in the answer $a$ corresponding to question $x$. We will simply write $S$ when there is no chance for confusion. We will write $S_x$ for the variables corresponding to question $x$. 
\end{definition}
One should think of this as an interactive proof system: There are two players Alice and Bob, and a verifier. The verifier samples a pair of questions $(x,y)$ according to $q$ from $Q \times Q$, then sends $x$ to Alice and $y$ to Bob. Alice and Bob then both send $m$ bit answers $a$ and $b$ to the verifier, respectively. The verifier evaluates $D(a,b|x,y)$ and Alice and Bob win when this is $1$.

We are using the notation of \cite{bowen_aldous--lyons_2024} here, specifically with regards to the variable set. A general game can have different question sets for Alice and Bob, and different answer sets for them too. However we will not need that level of generality.

\subsection{Strategies}
Note that the only thing Alice and Bob can do is when asked questions $(x,y)\in Q^2$, give answers $(a,b) \in (\{0,1\}^m)^2$. So a strategy for Alice and Bob will be some conditional distribution $(p(a,b|x,y))_{ (\{0,1\}^m)^2\times Q^2}$. We will call these \textbf{correlations}. The idea is in these games, Alice and Bob can discuss a strategy before hand but cannot communicate during the game. So Alice has no idea what question Bob got and vice versa. So strategies where Alice's answer is independent of Bob's answer would be a classical strategy. 

However, if Alice and Bob share quantum entanglement, then there can be correlation between their choices of their answer. Precisely:
\begin{definition}
    Let $\mathfrak{G}$ be a non-local game. A \textbf{quantum strategy} is a map $\rho: S(\mathfrak{G})\to U_n(\mathbb{C})$ for some $n$ with:
    \begin{enumerate}
        \item $\rho(u)^2=1$ for each $u\in S$
        \item $\rho(u)$ and $\rho(v)$ commute whenever $u,v\in S_x$, i.e different bits corresponding to the same answer should commute.
    \end{enumerate}
    Define
    \[e_x^a = \prod_{i\leq m} \dfrac{1 + (-1)^{a_i}\rho(u_{x,i})}{2}.\]
    This is the product of the spectral projections of each bit onto $(-1)^{a_i}$. Then the corresponding correlation is defined as:
    \[p(a,b|x,y) = \tr_n(e_x^a e_y^b).\]
    A \textbf{quantum approximate correlation} is a pointwise limit of quantum correlations.
    
    The value of a game $\mathfrak{G}$ on strategy $\rho$, $\operatorname{Val}(\mathfrak{G},\rho)$, is just the expected value of winning if Alice and Bob answer $(a,b)$ on questions $(x,y)$ with probability $p(a,b|x,y)$. We define the \textbf{quantum value} of the game as:
    \[\omega^*(\mathfrak{G}):= \sup_{\rho}\operatorname{Val}(\mathfrak{G},\rho)\]
    the best probability of winning with a quantum strategy.
\end{definition}
Actually in literature these are called synchronous quantum approximate strategies, but we will not consider any non-synchronous strategies so we omit that word. This is also defined slightly differently than standard literature, but Theorem \ref{thm:3.3} shows the equivalence. Note that if the dimension was $1$, this would just be a classical correlation, since we would have $p(a,b|x,y) = p(a|x)p(b|y)$.

%\mbnote{You should define strategies $p(a,b|x,y)$ as well, before giving this definition.  As this point $p(a,b|x,y)$ and expected value are  not clearly defined.  Also, it would be worthwhile adding a remark that one-dimensional q-strategies are deterministic strategies.}
%\mbnote{Should you say what value 1 means?  i.e., perfect strategy?}

There are two perspectives, one of them is we can associate to each bit of the answer a formal variable, to be interpreted as matrices. This is what we are doing. In this case, the spectral projections of these variables correspond to what bit that variable outputs. On the other hand, one can associate a variable $v_x$ to each question $x\in Q$ with $v_x^{2^m}=1$. So now $v_x$ will have $2^m$ spectral projections, and they will correspond to the $2^m$ different outputs on question $x$. These perspectives are equivalent:

Let $\mathfrak{G}$ be a non-local game with question set $Q$, answers of size $m$, and variable set $S$. By $\mathcal{F}(S,2)$ we denote the group generated by elements of $S$ with the condition they have order $2$. By $\Conv_Q \mathbb{Z}_2^m$ we will denote the quotient of $\mathcal{F}(S,2)$ by $\{[u_{x,i}, u_{x,j}]: x\in X, i,j\leq m\}$. That is, we are imposing the condition that the variables corresponding to the same question commute. Let $U_{x,i}$ denote the group unitaries corresponding to $u_{x,i}$ in $C^*(\Conv_Q \mathbb{Z}_2^m)$, and $e_x^a$ the spectral projections as above for $a\in \mathbb{Z}_2^m$. That is: 
\[e_x^a = \prod_{i\leq m} \dfrac{1 + (-1)^{a_i}U_{x,i}}{2}\]
Note that these projections generate the algebra and have no relation between each other except $e_x^a e_x^b = \delta_{ab}$ and $\sum_a e_x^a = 1$. Hence we can find an isomorphism $C^*(\Conv_Q \mathbb{Z}_2^m) \cong C^*(\mathcal{F}(Q, 2^m))$ by looking at the $2^m$ spectral projections of the defining unitaries in the latter. In particular, this algebra has the local lifting property (See \cite[Corollary 4.6]{enders_local_2024}). 

Any set of projections satisfying the same properties as $e_x^a$ are called \textbf{projective valued measurements (PVMs)}. One should think of these as ``quantum probabilities'' of answering $a$ on question $x$. To get actual probabilities out of these, one needs to use traces:

\begin{theorem}\cite[Theorem 3.6]{kim_synchronous_2018}\label{thm:3.3}
    Let $\mathfrak{G}$ be a non-local game with question set $Q$ and answers of size $m$. Then $p$ is a quantum approximate correlation iff there is an amenable trace $\tau$ on $C^*(\Conv_Q \mathbb{Z}_2^m) \cong C^*(\mathcal{F}(Q,2^m))$ so that
    \[p(a,b|x,y) = \tau(e_x^a e_y^b).\]
\end{theorem}
If we remove the restriction to amenable traces and allow any trace, then the resulting class of correlations are called quantum commuting. The idea of $\mathsf{MIP}^*=\mathsf{RE}$ was separating the quantum commuting value and quantum value of some game, and that would give a non-amenable trace on this group. 

Of course, we do not want just any old trace to induce our correlations, but specifically the ones coming from IRS's. So let us define:
\begin{definition}
    Let $\mathfrak{G}$ be a non-local game with question set $Q$ and answers of size $m$. An \textbf{IRS strategy} is defined by the following data: a map $\sigma:S\cup \{J\}\to \operatorname{Aut}(X,\mu)$ where $(X,\mu)$ is a Borel probability space. We will require (equality here means they agree $\mu$-almost everywhere):
    \begin{itemize}
        \item $\mu(x: \sigma(J)x=x)=0$ and $\sigma(J)^2=1$.
        \item For $u\in \sigma(S)$, we have $u^2=1$ and $u\sigma(J)=\sigma(J)u$.
        \item For $u,v \in \sigma(S_x)$ we have $uv=vu$.
    \end{itemize}
     If $\mu$ was the uniform measure on a finite set, call this a \textbf{permutation strategy}.
\end{definition}
Note that $\sigma$ extends to a homomorphism $\tilde{\sigma}:\Conv_Q \mathbb{Z}_2^m \times \mathbb{Z}_2 \to \Aut(X,\mu)$. That is, it is a p.m.p. action of the group $\mathbb{Z}_2^m \times \mathbb{Z}_2$ on $(X,\mu)$. So we can consider the IRS $\Stab(\tilde{\sigma})$ and then the corresponding trace $\tau$ on $C^*(\Conv_Q \mathbb{Z}_2^m \times \mathbb{Z}_2)\cong C^*(\mathcal{F}(Q,2^m))\otimes C^*(\mathbb{Z}_2)$ has 
    \[\tau(g) = \mu(x: gx=x).\] 
    %\mbnote{You should note that $\sigma$ induces an action $(g,x) \mapsto gx$ of this group on $X$, before the formula is presented.}
    So we can define the correlation:
    \[p(a,b|x,y) = \tau((1-J) e_x^ae_y^b).\]
    Note that this is positive as $\tau((1-J)a) = \frac{1}{2} \tau((1-J)a(1-J))$ and these sum up to $1$ as $\tau(J)=0$. We will define $\omega_{IRS}(\mathfrak{G})$ to be the best winning probability over all IRS strategies. 

%\mbnote{This is a superlong defn.  It should not include all the discussion - that should be part of the regular body of the paper.}
%\mbnote{Also a reminder about the $X = Q$ thing.  Actually it seems $Q$ is better for the question set, as now $X$ is the measure space.}  
We remark that this definition agrees with that of \cite[Definition 6.11]{bowen_aldous--lyons_2024-1} for permutation strategies. They defined permutation strategies as maps $\sigma: S\cup\{J\} \to \operatorname{Sym}(2n)$ with the same conditions, and if $W^{-}$ is the negative eigenspace of the image of $J$ in $U(2n)$, then they defined the correlation as 
\[p(a,b|x,y) = \tr_n (e_x^ae_y^b|_{W^{-}}).\]
Since $\frac{1-J}{2}$ is the projection onto $W^{-}$, this is the same correlation as our definition. This is why we had to add this $(1-J)$ to our construction, as without it we would not be able to use the result TailoredMIP*=RE of \cite{bowen_aldous--lyons_2024}. What we have done is in essence define the dual of the permutation strategies of \cite{bowen_aldous--lyons_2024-1}, which will allow us a direct separation of the convex sets in question without going through some intermediary like subgroup tests.

Now if we had a game with a IRS strategy that is not quantum approximate, then intuitively that should arise from an IRS trace that is not amenable, which is what we want. Formally:
\begin{theorem}\label{thm:3.5}
    Suppose there is a game $\mathfrak{G}$ with 
    \[\omega_{IRS}(\mathfrak{G})>\omega^*(\mathfrak{G}).\]
    Then there is a non co-hyperlinear IRS on some free group.
\end{theorem}
\begin{proof}
    Let $\mathfrak{G}$ be a game with $\omega_{IRS}(\mathfrak{G})>\omega^*(\mathfrak{G})$. Let $\sigma:S\cup \{J\}\to \Aut(X,\mu)$ be a IRS strategy with $\operatorname{Val}(\mathfrak{G}, \sigma)>\omega^*(\mathfrak{G}).$ Let $\tau$ be the corresponding trace on $C^*(\mathcal{F}(Q,2^m))\otimes C^*(\mathbb{Z}_2)$ so that
    \[p(a,b|x,y) = \tau((1-J)e_x^a e_y^b)\] is the correlation.

    Note that there is a natural inclusion $ C^*(\mathcal{F}(Q,2^m))\subset C^*(\mathcal{F}(Q,2^m))\otimes C^*(\mathbb{Z}_2)$ induced by the natural group inclusion\cite[Proposition 2.5.8]{brown_c-algebras_2008}. Note that $a\mapsto \tau((1-J)a)$ for $ a \in C^*(\mathcal{F}(Q,2^m))\otimes C^*(\mathbb{Z}_2)$ is a trace as $J$ commutes with everything. Call this trace $\tau'$. Then $\tau'|_{C^*(\mathcal{F}(Q,2^m))}$ is also a trace, and has
    \[p(a,b|x,y) = \tau'|_{C^*(\mathcal{F}(Q,2^m))} (e_x^a e_y^b).\]
    Since $p$ is not a quantum approximate correlation, $\tau'|_{C^*(\mathcal{F}(Q,2^m))}$ is not an amenable trace. Since amenability of traces is preserved under restriction \cite[proposition 6.3.5]{brown_c-algebras_2008}, we get $\tau'$ itself is not amenable. Now we have \[\tau(a) = \frac12 \tau((1-J)a) + \frac12 \tau(Ja)\] and both are traces, and since the amenable traces are a face \cite[proposition 6.3.7]{brown_c-algebras_2008}, we have $\tau$ is not amenable.

    Note that $C^*(\Conv_Q \mathbb{Z}_2^m \times \mathbb{Z}_2) \cong  C^*(\mathcal{F}(Q,2^m))\otimes C^*(\mathbb{Z}_2)$ has the local lifting property \cite[exercise 13.2]{brown_c-algebras_2008}. Note there is a natural quotient map $q: \mathcal{F}(S\cup\{J\}) \to \Conv_Q \mathbb{Z}_2^m \times \mathbb{Z}_2$ and it induces a quotient on the $C^*$ level. Here $\mathcal{F}(S\cup\{J\})$ is the free group generated by $S\cup \{J\}$, So we can lift up $\tau$ to $\tilde{\tau} = \tau\circ q$ on $C^*(\mathcal{F}(S\cup\{J\}))$. Note $\tilde{\tau}$ corresponds to an action of $\mathcal{F}(S\cup\{J\})$ on $(X,\mu)$, so it is an IRS trace. We also have 
    \[\tilde{\tau}(\mathcal{F}(S\cup\{J\}))'' = \tau\Bigl(\Conv_Q \mathbb{Z}_2^m \times \mathbb{Z}_2\Bigr)''.\]
    Since the latter has the local lifting property, the non-amenability of $\tau$ implies the algebra is not Connes embeddable, and hence that $\tilde{\tau}$ is not amenable. 
\end{proof}
Since this implies the negation to Aldous-Lyons by Proposition \ref{prop:2.22}, this proof allows us to skip  \cite{bowen_aldous--lyons_2024-1}. That is, we reduce to non-local games without having to go through subgroup tests. In that paper, they also reduced to the a priori harder problem of finding this separation for a restricted class of games called tailored games, which have a sort of linearization baked into them allowing translation as group like objects. However we do not require this, as these IRS strategies already bake in a group like structure.

This highlights two different approaches when one wants to use these non-local games: restrict the class of strategy or restrict the class of games. \cite{bowen_aldous--lyons_2024-1} ended up doing both but their ideas were based off of the latter. If there was a linear constraint game with perfect quantum commuting strategy but $\omega^*<1$, then there would be a non-hyperlinear group. Based off of this, the tailored games of \cite{bowen_aldous--lyons_2024-1} were designed to be ``semi-linear''. 

%\mbnote{Some of the discussion of the previous paragraphs comparing with BCLV approach should also appear in the introduction.}

\subsection{NPA hierarchy}\label{sec:3.2}
To get the separation, we will need an analogue of the NPA hierarchy for these IRS strategies. That is:
\begin{theorem}\label{thm:3.6}
    Let $\mathfrak{G}$ be a non-local game, then there is a computably enumerable sequence $\alpha_n$ that is monotonically decreasing and converges to $\omega_{IRS}(\mathfrak{G})$.
\end{theorem}
Recall this means there is a Turing machine which takes as input the description of some non local game, and enumerates the sequence $\alpha_n$.

To see this we will recall some notions from \cite[section 2]{bowen_aldous--lyons_2024-1}. Our NPA heierarchy for IRS value of a game is inspired by the one in \cite{bowen_aldous--lyons_2024-1} for subgroup tests, and is using some of their lemmas leading upto that.
Throughout this section, fix a non-local game $\mathfrak{G}$ and its variable set $S$. Then let $\mathcal{F}=\mathcal{F}(S\cup\{J\})$ the free group generated by $S\cup\{J\}$.
\begin{definition}
   Let $B\subset C\subset \mathcal{F}$. Then define the restriction map $R_{B\subset C}: \{0,1\}^C\to \{0,1\}^B$ as the map that sends $A\subset C$ to $A\cap B \subset B$. If $C=\mathcal{F}$, then call this map $R_B$.
\end{definition}
\begin{definition}
    A pseudo subgroup of $B\subset \mathcal{F}$ is a subset $A\subset B$ which is the restriction of some subgroup. That is, there is some subgroup $H\subset \mathcal{F}$ so that $A=R_B(H)= H\cap B$.
\end{definition}
\begin{definition}
    A \textbf{random pseudo-subgroup} of $B\subset \mathcal{F}$ is a probability measure $\pi$ on $\{0,1\}^B$ supported on the pseudo-subgroups of $B$.

    A random pseudo-subgroup $\pi$ is in addition invariant if for each $K,L\subset B$ so that $sKs^{-1}, sLs^{-1}\subset B$ for each $s\in S\cup S^{-1}$, we have
    \[\pi\{A\subset B: K\subset A,\, A\cap L =\varnothing\} = \pi\{A\subset B: sKs^{-1}\subset A,\, A\cap sLs^{-1} =\varnothing\}\]

    Denote the set of Invariant random pseudo-subgroups of $B$ as $\mathcal{Q}_B$, and let $\tilde{\mathcal{Q}}_B = R_{B*}^{-1} (\mathcal{Q}_B)$. That is, the set of measures on $\{0,1\}^\mathcal{F}$ whose pushforward to $\{0,1\}^B$ is an invariant random pseudo-subgroup.
\end{definition}
For each $\pi \in \tilde{\mathcal{Q}}_B$, we will extend the value of a non-local game to it. Here is the idea, if $\sigma: S\cup \{J\} \to \Aut(X,\mu)$ is an IRS strategy, then it corresponds to an IRS $H=\Stab(\sigma)$ on $\mathcal{F}$. Now the correlations with respect to this are
\[p(a,b|x,y) = \tau_H((1-J) e_x^ae_y^b), \quad e_x^a = \prod_{i\leq m} \dfrac{1+(-1)^{a_i}U_{x,i}}{2}.\]
In particular, each $p(a,b|x,y)$ only requires the data of $\tau$ on a finite number of group elements. One can see this by expanding out $e_x^a$ and $e_y^b$. The only group elements appearing in this will be those that are a product of elements $J, U_{x,i}, U_{y,i}$ of length $\leq m+1$.

So we can find a finite set $K$ with 
\[p(a,b|x,y) = \sum_{g\in K} c_g^{abxy} \mathbb{P}(g\in H).\]
We can assume $K$ contains $S\cup \{J\}$, all of the squares from $S\cup \{J\}$ and all the commutators from it, simply by unioning this finite set with $K$.

Now for a probability distribution $\pi$ on $\{0,1\}^B$ for some set $B\supset K$, we can define $\operatorname{Val}(\mathfrak{G},\pi)$ as the expected value of winning with correlation
\[p(a,b|x,y) = \sum_{g\in K} c_g^{abxy} \pi(A: g\in A).\]

Note this will typically not actually be a conditional probability distribution, as these can become negative. This is because $\pi$ is not actually invariant and so we are not dealing with traces. But still we can find the value of the game on it. That is if $q$ is the distribution on the question set $Q\times Q$ and $D$ is the decider, then define
\[\operatorname{Val}(\mathfrak{G},\pi) = \sum_{x,y}\sum_{a,b} q(x,y)D(a,b|x,y) \biggr( \sum_{g\in K} c_g^{abxy} \pi(A: g\in A).\biggl).\]

Now we prove:
\begin{lemma}
    Let $K\subset B\subset C\subset \mathcal{F}$, and let $\pi$ be a probability measure on $\{0,1\}^C$. Then 
    \[\operatorname{Val}(\mathfrak{G},\pi) = \operatorname{Val}(\mathfrak{G},R_{B\subset C*}\pi).\]
\end{lemma}
\begin{proof}
    Note that 
    \[R_{B\subset C*}(\pi) (A\subset B : g\in A)  = \pi(A \subset C : g \in R_{B\subset C}(A)) = \pi(A \subset C : g \in A)\]
    for $g\in K$. Note we used that $g\in B$ always. Since these determine the value of the game, we are done.
\end{proof}
Only one thing is missing: not all IRS will define a strategy for our game. For example we want all members of $S$ to square to identity in the action. But this is easy to enforce:
\begin{definition}
    Let $K\subset B\subset \mathcal{F}$. Define the set $\mathcal{T}_B$ to be the set of probability measures on $\{0,1\}^B$ that has measure $1$ on the following sets:
    \begin{enumerate}
        \item $\{A: J\not\in A\}$
        \item $\{A: [u,v] \in A \text{ for } u,v\in S_{x}\}$ for each $x\in Q$.
        \item $\{A: u^2\in A\text{ for } u\in S\cup\{J\}\}$
        \item $\{A: [u,J] \in A \text{ for } u\in S\}$.
    \end{enumerate}
    Note $R_{B\subset C*}^{-1} \mathcal{T}_B= \mathcal{T}_C$.
\end{definition}
The point of this is that $\IRS(\mathcal{F})\cap \mathcal{T}_\mathcal{F}$ correspond precisely to the IRS's which induce strategies on $\mathfrak{G}$. To see this take any IRS $H$ of $\mathcal{F}$, and any action $\sigma$ of $\mathcal{F}$ on some Borel probability space $(X,\mu)$ which has $H=\Stab(\alpha)$. Then having measure $1$ on e.g. $\{A: u^2\in A\text{ for } u\in S\cup\{J\}\}$ is simply saying $u^2$ stabalizes almost all points of $X$, i.e $u^2=1$.

Now we can finally construct the NPA hierarchy (the proof follows the structure of \cite[proof of theorem 1.10]{bowen_aldous--lyons_2024-1}, with the appropriate modification for our setting):
\begin{proof}[Proof of Theorem \ref{thm:3.6}]
    Let $K$ be as above and $B\supset K$ be a finite set. Then we can define
    \[\operatorname{Val}_B(\mathfrak{G}) = \sup(\operatorname{Val}(\mathfrak{G},\pi): \pi \in \mathcal{Q}_B\cap \mathcal{T}_B) =\sup(\operatorname{Val}(\mathfrak{G},\pi): \pi \in \tilde{\mathcal{Q}}_B\cap \mathcal{T}_\mathcal{F}).\]
    Here we use that $R_{B*}^{-1} (\mathcal{T}_B) =\mathcal{T}_\mathcal{F}$. Let $B_t$ be a sequence of sets, $B_1\supset K$ and $B_{t+1}\supset B_t$. Also suppose that $\bigcup_t B_t = \mathcal{F}$. Then we prove:
    \begin{itemize}
        \item $\Val_{B_t}(\mathfrak{G})$ is computable. To see this, first from \cite[lemma 2.16]{bowen_aldous--lyons_2024-1}, $\mathcal{Q}_{B_t}$ can be considered as a computable polytope of $\mathbb{R}^{\{0,1\}^{B_t}}$ (i.e assigning each subset of $B_t$ a real number). This means $\mathcal{Q}_{B_t}$ can be computably defined using integer linear equations and inequalities. Note that $\mathcal{Q}_{B_t}\cap \mathcal{T}_{B_t}$ just adds finitely more linear equations to this: for fixed subsets of $B_t$ we are asking the sum over subsets of those to add up to $1$. Hence we can do a linear program over $\mathcal{Q}_{B_t}\cap \mathcal{T}_{B_t}$ to compute $\Val_{B_t}(\mathfrak{G})$ in finite time.
        \item $\Val_{B_t}(\mathfrak{G})\geq \Val_{B_{t+1}}(\mathfrak{G})$. To see this, note that $R_{B\subset C*}(\mathcal{Q}_C) = \mathcal{Q}_B$. So 
        \[\tilde{\mathcal{Q}}_{B_t} \cap \mathcal{T}_\mathcal{F} \supset\tilde{\mathcal{Q}}_{B_{t+1}} \cap \mathcal{T}_\mathcal{F}\]
        giving the result
        \item $\Val_{B_t}(\mathfrak{G}) \to \omega_{IRS}(\mathfrak{G})$. The $\tilde{\mathcal{Q}}_B$ and $\mathcal{T}_\mathcal{F}$ are weak$^*$ closed subsets of $\operatorname{Prob}(\{0,1\}^\mathcal{F})$ and hence their intersection is compact. So there is some $\pi_t \in \tilde{\mathcal{Q}}_{B_t}\cap \mathcal{T}_\mathcal{F}$ which has
        \[\Val_{B_t}(\mathfrak{G}) = \Val(\mathfrak{G}, \pi_t).\]
        We can pick a weak$^*$ cluster of $\pi_t$, say $\pi_\infty$. By passing to a subsequence, we may assume $\pi_t$ weak$^*$ converge to $\pi_\infty$. Now note that for all $s\geq t$, we have $\pi_s \in \tilde{\mathcal{Q}}_{B_t}\cap \mathcal{T}_S$. Since this set is closed, we get $\pi_\infty \in \tilde{\mathcal{Q}}_{B_t}\cap \mathcal{T}_\mathcal{F}$ for all $t$. I.e
        \[\pi_\infty \in\bigcap_t  \tilde{\mathcal{Q}}_{B_t}\cap \mathcal{T}_\mathcal{F}= \IRS(\mathcal{F}) \cap \mathcal{T}_\mathcal{F}.\]
        See \cite[lemma 2.12]{bowen_aldous--lyons_2024-1} for the proof of this intersection.
        
        This means $\pi_\infty$ is a IRS strategy for the game. And now since $\Val(\mathfrak{G},-)$ is a continuous functional we get:
        \[\lim_{t\to \infty} \Val_{B_t}(\mathfrak{G}) = \Val (\mathfrak{G}, \pi_t) = \Val(\mathfrak{G},\pi_\infty )\leq \omega_{IRS}(\mathfrak{G}).\]
        But as $\tilde{\mathcal{Q}}_{B_t}\cap \mathcal{T}_\mathcal{F}\supset \IRS(\mathcal{F}) \cap \mathcal{T}_\mathcal{F}$, we have $\Val_{B_t}(\mathfrak{G})\geq \omega_{IRS}(\mathfrak{G})$. We are done.
    \end{itemize}
\end{proof}
\subsection{Separation}
We now get the separation directly from the main theorem of \cite{bowen_aldous--lyons_2024}:
\begin{theorem}\label{thm:3.14}
    For each Turing machine $\mathcal{M}$, there is a non-local game $\mathfrak{G}_\mathcal{M}$ whose description can be computed in finite time from $\mathcal{M}$ with:
    \begin{itemize}
        \item If $\mathcal{M}$ halts, then $\mathfrak{G}_\mathcal{M}$ has a perfect permutation strategy.
        \item If $\mathcal{M}$ does not halt, then $\omega^*(\mathfrak{G}_\mathcal{M})<\frac{1}{2}$.
    \end{itemize}
\end{theorem}
Actually they proved something a priori stronger, that this is true when restricting to a class of games called tailored non-local games. With this we can prove:
\begin{theorem}\label{thm:3.15}
    There exists a non local game $\mathfrak{G}$ with $\omega_{IRS}(\mathfrak{G})> \omega^*(\mathfrak{G})$.
\end{theorem}
\begin{proof}
    First note for each $\mathfrak{G}$, there is a computably enumerable sequence $\beta_n$ that is increasing and converges to $\omega^*(\mathfrak{G})$ .This is the inclusion $\mathsf{MIP}^*\subset \mathsf{RE}$ \cite[Corollary 12.9]{ji_mipre_2022}. This is by exhaustive search: a quantum strategy only requires the the choice of unitaries of order $2$ with some commuting properties, and we can enumerate the finite set of rational unitaries in dimension $n$ with denominators atmost $n$. The density of rational unitaries give the convergence \cite[Theorem 2]{LuoDengChenYang2013}.
    
    Suppose for all games we had $\omega^*(\mathfrak{G})\geq \omega_{IRS}(\mathfrak{G})$. Under this assumption we will construct an algorithm for the Halting problem.

    Let $\mathcal{M}$ be a Turing machine, and $\mathfrak{G}_\mathcal{M}$ the game from Theorem \ref{thm:3.14}. On this input, our algorithm will alternate enumerating the $\alpha_n$ from Theorem \ref{thm:3.6} and the $\beta_n$ for $\mathfrak{G}_\mathcal{M}$. The algorithm will accept on the input if ever $\beta_n\geq \frac{1}{2}$. It will reject if ever $\alpha_n <1$.

    There are two cases: First suppose $\mathcal{M}$ halts. This means there is a perfect permutation strategy and so $\omega_{IRS}(\mathfrak{G}_\mathcal{M})=1$ and $\omega^*(\mathfrak{G}_\mathcal{M})=1$. So $\alpha_n$ are always $\geq 1$ and the $\beta_n$ are eventually $\geq 1/2$. So the algorithm accepts on this input.

    Suppose $\mathcal{M}$ never halts. This means $\omega^*(\mathfrak{G}_\mathcal{M})<\frac{1}{2}$ and so $\omega_{IRS}(\mathfrak{G}_\mathcal{M})<\frac12$. In particular, $\alpha_n$ is eventually $<1$ and $\beta_n$ are never $\geq 1/2$. So the algorithm rejects on this input.
    
\end{proof}

We now combine the results of the paper to prove the main theorem:
\begin{proof}[Proof of Theorem \ref{thm:1.1} and \ref{thm:1.0}]
    By Theorem \ref{thm:3.5} and \ref{thm:3.15} we see that there is a non co-hyperlinear IRS $H$ on some free group $\Gamma$. 
    
    Since $\IRS(\Gamma)$ and the set of co-hyperlinear IRS are closed convex sets, a separation between them implies some extreme point of $\IRS(\Gamma)$ is not co-hyperlinear. So there is some ergodic $H\in \IRS(\Gamma)$ which is not co-hyperlinear.
    
    By Proposition \ref{prop:2.17}, this means $L(\Gamma/H)$ is not Connes embeddable. By the discussion proceeding Definition \ref{def:1.14}, we have $L(\Gamma/H)$ embeds into some $L(\mathcal{R})$ where $\mathcal{R}$ is an countable p.m.p. relation. $\mathcal{R}$ can be chosen to be ergodic as $H$ was ergodic.
    
    Finally by \cite[Theorem D]{chifan_embedding_2022} we have $L(\mathcal{R})$ embedds into some $L(\mathcal{R}')$, where $\mathcal{R}'$ is a ergodic countable p.m.p. relation with property $(T)$. Since $L(\Gamma/H)$ was not Connes embeddable, $L(\mathcal{R}')$ is not either. 
\end{proof}

%% file: further.tex
\section{Further Directions}
Here are some future research directions:
\begin{enumerate}
    \item Given a trace $\tau$ on a free group (or any group), is there any way to determine if $\tau$ arises from an IRS? An intrinsic characterization of this would be interesting. Maybe one could try to generalize the notion of an IRS trace to general $C^*$ algebras if such a characterization existed.

    \item \cite{bowen_aldous--lyons_2024-1} proved the existence of a non co-sofic IRS on free groups, and we proved the existence of a non co-hyperlinear IRS. So a natural question is, are these notions the same? If they are the same, then hyperlinearity and soficity is the same for groups.

    If they are distinct, a naive application of non-local games will not help. The non-local game method for seperating two convex sets $B\subset A$ requires one to be able to have computable upper bounds for $A$ and computable lower bounds for $B$. The issue is $\omega^*$ cannot have computable upper bounds, otherwise it would be computable. So really, something novel is needed for this. Such a separation would be evidence for soficity and hyperlinearity not agreeing on groups.

    It is not clear that they should be distinct either, perhaps ergodic theoretic methods could be used to show they agree?

    We note that the existence of a Hilbert Schmidt stable group that is not permutation stable \cite{eckhardt_amenable_2023} implies a non co-sofic co-hyperlinear IRS on that group, but this does not extend to free groups. 

    \item 
    Let $\mathcal{R}$ be a p.m.p equivalence relation on $(X,\mu)$. Consider the uniform metric $d_u$ on the full group $[\mathcal{R}]\subset \Aut(X,\mu)$. That is
    \[d_u(T,S) = \mu(x: T(x)\neq S(x)).\]
    It is not hard to see that $L(\mathcal{R})$ being Connes embeddable is the same as $[\mathcal{R}]$ being metricly hyperlinear. Compare with \cite[Proposition 2.4]{cordeiro_elementary_2016}.

    However the issue is a metric group can fail to be metrically hyperlinear while the underlying discrete group is hyperlinear. 
    Can we drop the metric condition in the above? If so, we have a non-hyperlinear group. 
    
    We expect this should not be possible: there should be a non-Connes embeddable $L(\mathcal{R})$ with hyperlinear full group. However, it is not clear how to prove it one way or the other.
\end{enumerate}

%% file: acknowledgements.tex
\section*{Acknowledgements}
    I would like to thank Michael Brannan for discussions, supervision and very thoroughly proofreading this paper. I would also like to thank Michael Chapman for discussion, proofreading and providing suggestions to significantly clean up proof of Theorem \ref{thm:3.15}. I would like to thank Stefan Frunza for pointing out the explicit GNS construction in Proposition \ref{prop:1.13}. I would also like to thank Jesse Peterson for pointing out the proof of Lemma \ref{lemma:2.21}. I would also like to thank the anonymous referee for helpful comments that cleaned up the exposition in parts of this paper.